\documentclass{gtart}

\def\ifplaintex{\expandafter\ifx\csname documentclass\endcsname\relax}

\ifplaintex 
\else
\fi

\expandafter\ifx\csname beginpicture\endcsname\relax
\expandafter\ifx\csname documentclass\endcsname\relax
\input pictex \else
\input prepictex \input pictex \input postpictex \fi\fi

\def\gt{{\mathsurround=0pt\it $\cal G\mskip-2mu$eometry \&\ 
$\cal T\!\!$opology}}        

\def\gtp{{\mathsurround=0pt\it $\cal G\mskip-2mu$eometry \&\ 
$\cal T\!\!$opology $\cal P\!$ublications}}  

\def\lognumber#1{\def\thelognumber{#1}}
\def\volumenumber#1{\def\thevolumenumber{#1}}
\def\papernumber#1{\def\thepapernumber{#1}}
\def\volumeyear#1{\def\thevolumeyear{#1}}

\def\pagenumbers#1#2{\def\startpage{#1}\def\finishpage{#2}}
\def\published#1{\def\publishdate{#1}}
\def\proposed#1{\def\theproposer{#1}}
\def\seconded#1{\def\theseconders{#1}}
\def\received#1{\def\receiveddate{#1}}
\def\revised#1{\def\reviseddate{#1}}
\def\accepted#1{\def\accepteddate{#1}}

\def\asciiaddress#1{\def\theasciiaddress{#1}}

\long\def\asciiabstract#1{\long\def\theasciiabstract{#1}}

\let\\\par\let\thelognumber\relax
\let\thevolumenumber\relax\let\thepapernumber\relax
\let\thevolumeyear\relax\let\thesamplenumber\relax\let\startpage\relax
\let\finishpage\relax\let\publishdate\relax\let\receiveddate\relax
\let\reviseddate\relax\let\accepteddate\relax\let\theasciititle\relax
\let\theasciiauthors\relax\let\theasciiaddress\relax
\let\theasciiabstract\relax
\let\theasciiemail\relax\let\theshortauthors\relax\let\theshorttitle\relax

\long\def\maketitlep{   

\count0=\startpage

\gt\hfill      
\beginpicture
\setcoordinatesystem units <0.33truein, 0.33truein> point at 2.2 0.9
\setplotsymbol ({$\cal G$})
\plotsymbolspacing=9truept
\circulararc 315 degrees from 0 1 center at 0 0
\setplotsymbol ({$\cal T$})
\circulararc 315 degrees from 1 -1 center at 1 0
\endpicture
\break
{\small\ifx\thesamplenumber\relax 
Volume \else Sample
\fi\thevolumenumber\ (\thevolumeyear)
\startpage--\finishpage\nl
Published: \publishdate}
\vglue 0.5truein plus 0.4fil minus 0.1truein

{\parskip=0pt\leftskip 0pt plus 1fil\def\\{\par\smallskip}{\ifplaintex\large
\else\Large\fi\bf\thetitle}\par\medskip}   

\vglue 0pt plus 0.1fil 

{\parskip=0pt\leftskip 0pt plus 1fil\def\\{\par}{\sc\theauthors}
\par\medskip}

\vglue 0pt plus 0.1fil 

{\small\parskip=0pt\let\newline\\
{\leftskip 0pt plus 1fil\def\\{\par}{\sl\theaddress}\par}
\expandafter\ifx\theemail\relax    
\relax\else\vglue 5pt plus 0.02fil minus 2pt\def\\{\stdspace{\rm 
and}\stdspace} 
\cl{Email:\stdspace\tt\theemail}\fi
\ifx\theurl\relax                  
\relax\else\vglue 5pt plus 0.02fil minus 2pt\def\\{\stdspace{\rm 
and}\stdspace}
\cl{URL:\stdspace\tt\theurl}\fi\par}

\vglue 7pt plus 0.3fil minus 3pt

{\bf Abstract}
\vglue 5pt plus 0.1fil minus 2pt

\theabstract

\vglue 7pt plus 0.3fil minus 3pt

{\bf AMS Classification numbers}\quad Primary:\quad \theprimaryclass

Secondary:\quad \thesecondaryclass

\vglue 5pt plus 0.3fil minus 2pt

{\bf Keywords}\quad \thekeywords

\vglue 10pt plus 0.5fil minus 5pt

{\small  Proposed: \theproposer\hfill Received: \receiveddate\nl
Seconded: \theseconders\hfill 
\ifx\reviseddate\relax                         
Accepted: \accepteddate                        
\else
Revised: \reviseddate                          
\fi}
\eject
}       

\let\maketitlepage\maketitlep
\let\maketitle\maketitlepage

\font\phead=cmsl9 scaled 950
\font=cmsl9 scaled 1050
\font\pnum=cmbx10 scaled 913
\font\lnum=cmbx10 
\font\pfoot=cmsl9 scaled 950
\font=cmsl9 scaled 1050
\ifplaintex
\headline{\vbox to 0pt{\vskip -4.5mm\line{\small\phead\ifnum
\count0=\startpage ISSN 1364-0380 (on line)
1465-3060 (printed) \hfill {\pnum\folio}\else\ifodd\count0\def\\{ }%
\ifx\theshorttitle\relax\thetitle\else\theshorttitle\fi\hfill{\pnum\folio}
\else\def\\{ and }{\pnum\folio}\hfill\ifx\theshortauthors\relax\theauthors
\else\theshortauthors\fi\fi\fi}\vss}}
\footline{\vbox to 0pt{\vglue 0mm\line{\small\pfoot\ifnum\count0=\startpage
\copyright\ \gtp\hfill\else
\gt, Volume \thevolumenumber\ (\thevolumeyear)\hfill\fi}\vss
}}
\else
\makeatletter
\def\@oddhead{{\small\lhead\ifnum\count0=\startpage ISSN 1364-0380 (on line)
1465-3060 (printed) \hfill {\lnum\number\count0}\else\ifodd\count0
\def\\{ }\ifx\theshorttitle\relax \thetitle \else\theshorttitle\fi\hfill
{\lnum\number\count0}\else\def\\{ and }{\lnum\number\count0}
\hfill\ifx\theshortauthors\relax 
\theauthors\else\theshortauthors\fi\fi\fi}}\def\@evenhead{\@oddhead}
\def\@oddfoot{\small\lfoot\ifnum\count0=\startpage\copyright\ \gtp\hfill\else
\gt, Volume \thevolumenumber\ (\thevolumeyear)\hfill\fi}
\def\@evenfoot{\@oddfoot}
\makeatother
\fi

\newwrite\gtoutfile
\long\gdef\makeheadfile{  
{\def\\{, }\def\s{ }
\immediate\openout\gtoutfile head.xxx
\immediate\write\gtoutfile{To: math@arxiv.org}
\immediate\write\gtoutfile{Subject: put or rep NNNNN:pppp}
\immediate\write\gtoutfile{--text follows this line--}
\immediate\write\gtoutfile{Proxy-for: \ifx\theasciiauthors\relax
\theauthors\else\theasciiauthors\fi\s<\ifx\theasciiemail\relax\theemail\else\theasciiemail\fi>}
\immediate\write\gtoutfile{\noexpand\\}
\immediate\write\gtoutfile{Authors: \ifx\theasciiauthors\relax
\theauthors\else\theasciiauthors\fi}
{\def\\{ }\immediate\write\gtoutfile{Title: \ifx\theasciititle\relax
\thetitle\else\theasciititle\fi}}
\immediate\write\gtoutfile{Subj-class: GT or SG or MG etc}
\immediate\write\gtoutfile{MSC-class: \theprimaryclass\ifx\thesecondaryclass\relax\else, \thesecondaryclass\fi}
\immediate\write\gtoutfile{Journal-ref: Geom. Topol. \thevolumenumber
(\thevolumeyear) \startpage-\finishpage}
\immediate\write\gtoutfile{Comments: Published by Geometry and Topology at}
\immediate\write\gtoutfile{\s\s http://www.maths.warwick.ac.uk/gt/GTVol\thevolumenumber/paper\thepapernumber.abs.html}
\immediate\write\gtoutfile{\noexpand\\}
\immediate\write\gtoutfile{}
\ifx\theasciiabstract\relax
\immediate\write\gtoutfile{\theabstract}\else
\immediate\write\gtoutfile{\theasciiabstract}\fi
\immediate\write\gtoutfile{}
\immediate\write\gtoutfile{\noexpand\\}
\immediate\write\gtoutfile{}
\immediate\closeout\gtoutfile}}  

\def\maketitlepage{\maketitlep\makeheadfile}
\let\maketitle\maketitlepage

\lognumber{304}
\volumenumber{7}\papernumber{10}
\volumeyear{2003}
\pagenumbers{321}{328}
\received{24 January 2003}
\revised{11 April 2003}
\accepted{14 May 2003}
\published{19 May 2003}
\proposed{Walter Neumann}
\seconded{Cameron Gordon, David Gabai}

\usepackage{amsmath,amssymb,graphicx}

\def\color[#1]#2{\relax}

\newtheorem*{rigidity}{Rigidity Theorem}
\newtheorem{thm}{Theorem}[section]

\newtheorem*{thm*}{Theorem}

\newtheorem*{cor*}{Corollary}

\newtheorem*{prop*}{Proposition} 
 
\newtheorem{lem}[thm]{Lemma} 
\newtheorem*{lem*}{Lemma}

\newtheorem*{claim*}{Claim} 
 
\newtheorem*{fact*}{Fact}

\theoremstyle{definition}
 
\newtheorem*{dfn*}{Definition}

\theoremstyle{remark}
\newtheorem*{rem*}{Remark}
\newtheorem*{example*}{Example}

\newcommand{\rond}[1]{\overset{\scriptscriptstyle \circ}{#1}} 
\newcommand{\es}{\emptyset}
\renewcommand{\phi}{\varphi} 
\newcommand{\m} {^{-1}}

\newcommand {\ra} {\rightarrow}

\newcommand{\ie} {ie, }

\newcommand{\bbR} {{\mathbb{R}}}

\newcommand{\Aut} {\mathop{\mathrm{Aut}}}
\newcommand{\Fix}{\mathop{\mathrm{Fix}}}

\begin{document}
\title{A very short proof of Forester's rigidity result}
\author{Vincent Guirardel}

\address{Laboratoire E. Picard, UMR 5580, B\^atiment 1R2\\Universit\'e Paul 
Sabatier, 118 rte de Narbonne\\31062 Toulouse cedex 4, France}

\asciiaddress{Laboratoire E. Picard, UMR 5580, Batiment 1R2\\Universite Paul 
Sabatier, 118 rte de Narbonne\\31062 Toulouse cedex 4, France}

\email{guirardel@picard.ups-tlse.fr}

\begin{abstract}
The deformation space of a simplicial $G$--tree $T$ is the set of $G$--trees which can be obtained from $T$
by some collapse and expansion moves, or equivalently, which have the same elliptic subgroups as $T$.
We give a short proof of a rigidity result by Forester
which gives a sufficient condition for a deformation space  
to contain an $\Aut(G)$--invariant $G$--tree. 
This gives a sufficient condition for a JSJ splitting to be invariant under automorphisms of $G$.
More precisely, the theorem claims that a deformation space contains at most one 
strongly slide-free $G$--tree, where strongly slide-free means the following: 
whenever two edges $e_1,e_2$ incident on a same vertex $v$ are such that 
$G_{e_1}\subset G_{e_2}$, then $e_1$ and $e_2$ are in the same orbit under $G_v$.
\end{abstract}

\asciiabstract{The deformation space of a simplicial G-tree T is the
set of G-trees which can be obtained from T by some collapse and
expansion moves, or equivalently, which have the same elliptic
subgroups as T.  We give a short proof of a rigidity result by
Forester which gives a sufficient condition for a deformation space to
contain an Aut(G)-invariant G-tree.  This gives a sufficient condition
for a JSJ splitting to be invariant under automorphisms of G.  More
precisely, the theorem claims that a deformation space contains at
most one strongly slide-free G-tree, where strongly slide-free means
the following: whenever two edges e_1, e_2 incident on a same vertex v
are such that G_{e_1} is a subset of G_{e_2}, then e_1 and e_2 are in
the same orbit under G_v.}

\primaryclass{20E08}                
\secondaryclass{57M07, 20F65}              
\keywords{Tree, graph of groups, folding, group of automorphisms}

\maketitlepage

In \cite{For_deformation}, Forester introduced the notion of \emph{deformation}
for simplicial trees with a cocompact action of a group $G$, or equivalently, for splittings
of $G$ as a finite graph of groups.
A deformation consists in a sequence of \emph{collapse} and \emph{expansion} moves in the following sense:
a \emph{collapse move} consists in replacing an edge in a graph of groups corresponding to
an amalgamated product $A*_C C$ by a vertex with vertex group $A$, and an \emph{expansion move}
is the inverse operation.

Remember that a subgroup of $G$ is \emph{elliptic} in a $G$--tree $T$ if it fixes a point in $T$.
Forester proves that two cocompact simplicial $G$--trees can be deformed into one another if and only if
they have the same elliptic subgroups \cite[Theorem 1.1]{For_deformation}.
In terms of the geometric realization of the trees, this can also be reformulated by saying that
two $G$--trees $T$ and $T'$ can be deformed into one another if
and only if there is an equivariant  continuous  map from $T$ to $T'$ and one from 
 $T'$ to $T$.

This notion of deformation is interesting because the various JSJ splittings introduced by Rips--Sela, Dunwoody--Sageev,
Fujiwara--Papasoglu \cite{RiSe_JSJ,DuSa_JSJ,FuPa_JSJ} are unique up to deformation \cite{For_uniqueness}.
On the other hand, the JSJ splittings introduced by Bowditch for one-ended hyperbolic groups
and by Scott--Swarup for finitely presented groups are really unique, up to $G$--equivariant isomorphism of trees
 \cite{Bo_cut,SS_regular}. 
In particular, $\Aut(G)$ acts naturally by isometries on the corresponding simplicial tree,
or equivalently, any outer automorphism of $G$ is induced by an automorphism of the corresponding graph of groups.
Therefore, this allows one to understand the automorphism group of $G$ by understanding the automorphisms
of the JSJ splitting (see \cite{BaJi_automorphism,Lev_automorphisms}).

Forester's rigidity theorem gives a sufficient condition for the
existence of a \emph{canonical} point in a deformation space and
hence gives a criterion for a JSJ splitting \`a la Rips--Sela,
Dunwoody--Sageev or Fujiwara--Papasoglu to be invariant under
$\Aut(G)$.

In the sequel, we assume that all actions are without inversions (\ie no element exchanges the two endpoints of an edge)
since one can get rid of inversions by taking the first barycentric subdivision of $T$.
We will also assume that the actions are \emph{minimal}, \ie with no proper invariant subtree.
Note that if $T$ is not assumed to be minimal, but if at least one element of $G$ is not elliptic,
then $T$ contains a unique minimal invariant subtree.
We denote by $G_v$ (resp.\ by $G_e$)
the stabilizer of a vertex $v$ (resp.\ of an edge $e$).

\begin{dfn*}[{\rm(}Strongly slide-free, reduced $G$--tree\/{\rm)}]
  A $G$--tree is \emph{strongly slide-free} if it satifies the following condition:
if two edges $e_1,e_2$ having a common vertex $v$
          are such that $G_{e_1}\subset G_{e_2}$, then $e_1$ and $e_2$ are in the same orbit under $G_v$.

A $G$--tree $T$ is \emph{reduced} if one cannot perform a collapse on $T$, \ie if 
 for each edge $e$ incident on some vertex $v$ such that $G_e=G_v$, then the two endpoints
of $e$ are in the same orbit.    
\end{dfn*}

Note that in a minimal strongly slide-free $G$--tree, no vertex stabilizer can fix an edge, thus
a minimal strongly slide-free splitting is itself reduced.
In the following result, the \emph{deformation space} of a $G$--tree $T$, is the set of all
$G$--trees $T'$ which can be deformed into $T$.

\begin{rigidity}[Forester {\cite[Corollary 1.3]{For_deformation}}]
There is at most one strongly slide-free minimal $G$--tree in each deformation space.

More precisely, let $T,T'$ be two minimal simplicial $G$--trees which have the same elliptic subgroups.
Assume that $T$ is strongly slide-free and that $T'$ is reduced.
Then there is a $G$--equivariant isomorphism between $T$ and $T'$ (and the isomorphism is unique).
\end{rigidity}

\begin{cor*}[Forester \cite{For_uniqueness}]
  If a group $G$ has a JSJ splitting which is strongly slide-free, then this splitting is $\Aut(G)$--invariant.
\end{cor*}

In particular, the action of $\Aut(G)$ on the Bass--Serre tree of the JSJ splitting provides a splitting
of $\Aut(G)$ as a graph of groups.

This result extends an earlier result by Gilbert--Howie--Metaftsis--Raptis and Pettet essentially claiming that
a deformation space contains at most one
minimal strongly slide-free $G$--tree satisfying the additional assumption that if two adjacent edges have nested stabilizers,
then these stabilizers coincide \cite{GHMR_tree,Pettet_automorphism}.
This result is in turn an extension of a result by
Karrass--Pietrowski--Solitar applying to amalgamated products \cite{KPS_automorphisms}.

 A similar result in a different situation is proved in \cite{GL_automorphismes}:
 it is shown that in each deformation space, if there is a $G$--tree with cyclic edge stabilizers
 which is acylindrical,  then the deformation space contains a $2$--acylindrical $G$--tree with cyclic edge stabilizers, 
 and the set of such $2$--acylindrical $G$--trees is a simplex. This gives a way to produce an
$\Aut(G)$--invariant JSJ--splitting for torsion free commutative transitive groups.

The proof given in \cite{For_deformation} is quite long and involved. 
The goal of this note is to give a very short alternative proof, in the spirit of the proof in \cite{GHMR_tree}.

\section{Definitions}

We recall shortly a few definitions and elementary properties. Consider a $G$--tree $T$. 
Given a vertex $v\in V(T)$ and an edge $e\in E(T)$, we will denote by $G_v$ and $G_e$ their stabilizer.
If an element $\gamma$ has a fix point in $T$, $\gamma$ is called \emph{elliptic}, and $\gamma$ is called \emph{hyperbolic} otherwise.
Similarly, we say that a subgroup $H<G$ is \emph{elliptic} if it fixes a point in $T$.
Given an elliptic element $\gamma\in G$, the fix set $\Fix \gamma$ of $\gamma$ is a subtree of $T$ (and the same of course holds for a subgroup). 
Serre's Lemma claims that
if $\Fix \gamma\cap\Fix\gamma'=\es$, then $\gamma\gamma'$ is hyperbolic \cite[Corollary 1, section 6.5]{Serre_arbres}.

Given two disjoint subtrees $A,B\subset T$, the \emph{bridge} between $A$ and $B$ is the smallest arc joining $A$ to $B$:
it is the arc $[a,b]$ such that $a\in A$, $b\in B$, and
any arc joining a point of $A$ to a point of $B$ contains $[a,b]$. The \emph{projection} $p(x)$ of $x$ on $A$ is
the closest point to $x$ in $A$; if $x\notin A$, $[p(x),x]$ is the bridge between $A$ and $\{x\}$.

We will often blur the distinction between $T$ and its geometric realization, thus identifying the edge $e$
with endpoints $a,b$ to an homeomorphic copy $[a,b]$ of the interval $[0,1]$ in $\bbR$, 
while $\rond{e}$ will represent the open segment $(a,b)=[a,b]\setminus\{a,b\}$.
Note that if an element (or a subgroup) of $G$ fixes a point in $T$, then it fixes a vertex of $T$ 
and the fixed subtree of an element (or a subgroup) is a simplicial subtree (this uses the absence of inversion).

\section{Proof of Forester's rigidity Theorem}

\begin{proof}[Proof of the rigidity Theorem]
In a minimal strongly slide-free $G$--tree, no vertex stabilizer can fix an edge.
In particular, vertex stabilizers of $T$ fix no more than one vertex; thus vertex stabilizers 
are characterized as maximal elliptic subgroups of $G$.

Let's now define a $G$--equivariant map $f:T\ra T'$.
For each vertex $v\in V(T)$, choose equivariantly a vertex $f(v)\in V(T')$ fixed by $G_v$,
and extend $f$ linearly and equivariantly on edges.
First, the restriction of $f$ to $V(T)$ is injective: if $f(u)=f(v)$, then $\langle G_u,G_v\rangle$ is elliptic
in $T'$, hence it is also elliptic in $T$. Since vertex stabilizers of $T$ are maximal elliptic, one gets $G_u=G_v=\langle G_u,G_v\rangle$, so $u=v$.
Note that this implies that the image of every edge of $T$ is a non-degenerate arc in $T'$.

We will prove that $f$ is an isomorphism. 
Since $T'$ minimal, $f$ is onto (as a topological map: some vertices of $T'$
may have no preimage in $V(T)$).

The strongly-slide free condition gives the following fact (see \cite{GHMR_tree}):
\begin{lem}\label{lem_orbit}
Assume that $e_1,e_2\in E(T)$ are two edges sharing a common vertex $v$
and that $f(e_1)\cap f(e_2)$ is not reduced to one point.
Then $e_1$ and $e_2$ are in the same $G_v$--orbit and 
 $f(e_1)\cap f(e_2)$ is strictly contained in $f(e_1)$ (resp. in $f(e_2)$). 
\end{lem}

\begin{rem*}
  The lemma implies that $f(e_1)$ and $f(e_2)$ have the same length.
\end{rem*}

\begin{proof}
Consider the group $H=\langle G_{e_1},G_{e_2} \rangle< G_v$.

First assume that $H$ fixes only $v$ and argue towards a contradiction. Consider the vertex $w'$ at distance 1 from $f(v)$ on $f(e_1)\cap f(e_2)$.
Since $G_{w'}$ fixes a vertex in $T$, and since $H\subset G_{w'}$, $G_{w'}$ fixes $v$ (and only $v$).
Therefore, $G_{w'}\subset G_v\subset G_{f(v)}$. Since $T'$ is reduced, $f(v)$ and $w'$ are in the same orbit,
hence $w'$ has a preimage $w$ in the orbit of $v$ ($w$ is thus a vertex of $T$). 
Now $G_w\subset G_{w'}\subset G_v$, hence $w=v$ by maximality of vertex stabilizers, contradicting $f(w)=w'\neq f(v)$.

Thus $H$ fixes a vertex different from $v$. Since $H$ also fixes $v$, $H$ fixes an edge $e_3$ incident on $v$.
Since $G_{e_1},G_{e_2}\subset G_{e_3}$, the strongly-slide free condition says that $e_1,e_2,e_3$
are in the same $G_v$--orbit. 

Finally, if one had $f(e_1)\cap f(e_2)=f(e_1)$, then 
one would have $f(e_1)=f(e_2)$ since those two arcs have the same length (they are in the same orbit),
and $f$ would identify two vertices, a contradiction.
\end{proof}

\begin{lem}\label{lem_tripode}
Assume that $e_1,e_2,e_3$ are three consecutive edges in $T$, 
then $f(e_1)\cap f(e_2)\cap f(e_3)=\es$.   
\end{lem}

\begin{proof}
Denote by $v_1$ the common vertex of $e_1$ and $e_2$, and by $v_2$ the common vertex of $e_2$ and $e_3$ (note that $v_2\neq v_1$). 
Assume that $f(e_1)\cap f(e_2)\cap f(e_3)$ contains a point $p'$.
Then $f(e_1)$ must meet $f(e_2)$ in more than one point since otherwise,
$f(e_2)$ would be contained in $f(e_3)$, a contradiction.
By the previous lemma, there is an element $\gamma_1\in G_{v_1}$ sending $e_1$ on $e_2$,
so $\gamma_1$ fixes pointwise $f(e_1)\cap f(e_2)$, hence $\gamma_1$ fixes $p'$.
Similarly, there is an element $\gamma_2\in G_{v_2}$ sending $e_2$ on $e_3$,
and which fixes $p'$.

Now $\gamma_1\gamma_2$ is elliptic in $T'$ (it fixes $p'$) and is hyperbolic
in $T$ by Serre's Lemma since $\Fix \gamma_1\cap\Fix \gamma_2$ is empty in $T$
because neither $\gamma_1$ nor $\gamma_2$ fix $e_2$. 
This is a contradiction.
\end{proof}

The following lemma will say that the image under $f$ of a non-backtracking path $v_0,\dots,v_n$ cannot backtrack too much.
\begin{lem}[Backtracking lemma]
Consider a sequence of vertices $u_0,\dots,u_n$ in a tree $T$ such that
  \begin{enumerate}
  \item $u_{i}\neq u_{i+1}$ for all $i\in\{0,\dots,n-1\}$;\label{a}
  \item  $[u_{i-1}, u_i]\cap [u_i, u_{i+1}]$ is strictly contained\label{b}
in $[u_{i-1}, u_i]$ and in $[u_i, u_{i+1}]$ for each $i\in\{1,\dots,n-1\}$;
  \item $[u_{i-1}, u_i]\cap [u_i, u_{i+1}]\cap [u_{i+1}, u_{i+2}]=\es$ for each $i\in\{1,\dots,n-2\}$.\label{c}
  \end{enumerate}

Then for $|j-i|\geq 2$, $[u_{i-1},u_i]\cap[u_{j-1},u_{j}]=\es$.
\end{lem}

\begin{figure}[ht!]\small
\centering
\begin{picture}(0,0)%
\includegraphics{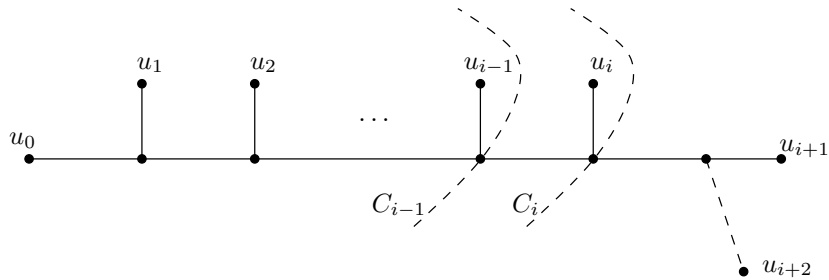}%
\end{picture}%
\setlength{\unitlength}{4144sp}%
\def\SetFigFont#1#2#3#4#5{\relax}
\begin{picture}(4680,1636)(136,-1010)
\put(901,254){\makebox(0,0)[b]{\smash{\SetFigFont{10}{12.0}{\rmdefault}{\mddefault}{\updefault}{\color[rgb]{0,0,0}$u_1$}%
}}}
\put(1576,254){\makebox(0,0)[b]{\smash{\SetFigFont{10}{12.0}{\rmdefault}{\mddefault}{\updefault}{\color[rgb]{0,0,0}$u_2$}%
}}}
\put(3601,254){\makebox(0,0)[b]{\smash{\SetFigFont{10}{12.0}{\rmdefault}{\mddefault}{\updefault}{\color[rgb]{0,0,0}$u_{i}$}%
}}}
\put(2926,254){\makebox(0,0)[b]{\smash{\SetFigFont{10}{12.0}{\rmdefault}{\mddefault}{\updefault}{\color[rgb]{0,0,0}$u_{i-1}$}%
}}}
\put(2251,-61){\makebox(0,0)[b]{\smash{\SetFigFont{12}{14.4}{\rmdefault}{\mddefault}{\updefault}{\color[rgb]{0,0,0}$\dots$}%
}}}
\put(4546,-961){\makebox(0,0)[lb]{\smash{\SetFigFont{10}{12.0}{\rmdefault}{\mddefault}{\updefault}{\color[rgb]{0,0,0}$u_{i+2}$}%
}}}
\put(4816,-241){\makebox(0,0)[b]{\smash{\SetFigFont{10}{12.0}{\rmdefault}{\mddefault}{\updefault}{\color[rgb]{0,0,0}$u_{i+1}$}%
}}}
\put(136,-196){\makebox(0,0)[b]{\smash{\SetFigFont{10}{12.0}{\rmdefault}{\mddefault}{\updefault}{\color[rgb]{0,0,0}$u_0$}%
}}}
\put(2566,-601){\makebox(0,0)[rb]{\smash{\SetFigFont{10}{12.0}{\rmdefault}{\mddefault}{\updefault}{\color[rgb]{0,0,0}$C_{i-1}$}%
}}}
\put(3241,-601){\makebox(0,0)[rb]{\smash{\SetFigFont{10}{12.0}{\rmdefault}{\mddefault}{\updefault}{\color[rgb]{0,0,0}$C_{i}$}%
}}}
\end{picture}
\caption{Backtracking lemma}
\label{fig}
\end{figure}

\begin{proof}
Let $C_i$ be the convex hull of $\{u_0,\dots,u_i\}$. We prove by induction the following property for $1\leq i\leq n-1$:
$$(P_i):\quad u_{i+1}\notin C_{i}\text{ and }[u_i,u_{i+1}]\cap C_{i-1}=\es.$$
The lemma will then follow immediately.

Since the property clearly holds for $i=1$, we prove $P_i\Rightarrow P_{i+1}$.
Assume that $u_{i+2}\in C_{i+1}$. Since $u_{i+2}\notin [u_{i},u_{i+1}]$ by hypothesis (\ref{b}),
$u_{i+2}$ lies in $C_{i-1}$ or in the bridge joining $[u_i,u_{i+1}]$ to $C_{i-1}$.
Thus $[u_{i+2},u_{i+1}]$ must meet the bridge between $[u_{i},u_{i+1}]$ and $C_{i-1}$,
hence must meet $[u_i,u_{i-1}]$, which contradicts hypothesis (\ref{c}).

If $[u_{i+1},u_{i+2}]$ meets $C_i$, then $[u_{i+1},u_{i+2}]$ contains the projection $p$ of $u_{i+1}$ on $C_i$.
Note that by definition $p\in[u_i,u_{i+1}]$, so $p\in[u_i,u_{i+1}]\cap[u_{i+1},u_{i+2}]$.
By (\ref{c}), $p\notin [u_{i-1},u_i]$. This implies that 
$p\in C_{i-1}$ since $p$ belongs to $C_i$ but not to the bridge joining $u_i$ to $C_{i-1}$.
Hence $p\in C_{i-1}\cap [u_i,u_{i+1}]$, which contradicts the induction hypothesis.
\end{proof}

Now let's conclude the proof of the theorem. 
Assume that $f$ is not an isomorphism. Then there exist two edges $e_1,e_2$ incident on a common vertex $v$ such that
$f(e_1)\cap f(e_2)$ contains more than one point.
Denote $v'=f(v)$, and let $w'\neq v'$ be the vertex at distance 1 from $v'$ on $f(e_1)\cap f(e_2)$.

Let $w$ be a point of $T$ (vertex or not) such that $f(w)=w'$. Denote by $v_0=v,v_1,\dots, v_n$ the vertices on the smallest simplicial
arc containing $[v,w]$ (in particular, $w\in [v_{n-1},v_n]$, and $w=v_n$ if and only if $w$ is a vertex).
Up to exchanging the roles of $e_1$ and $e_2$, we may assume that $[v_0,v_n]$ meets $e_1$ only at $v$.
Define $v_{-1}$ so that $[v_{-1},v_0]=e_1$. Thus, the vertices of the arc $[v_{-1},v_n]$ are $v_{-1},v_0,v_1,\dots, v_n$.
Lemma \ref{lem_orbit} and \ref{lem_tripode} say that one can apply 
the backtracking lemma to the sequence $u_i=f(v_{i})$, $i\in\{-1,0,\dots,n\}$.
Since $w'\in [u_{-1},u_0]\cap [u_{n-1},u_n]$, one gets $n=1$. By the Lemma \ref{lem_orbit}, the edge $[v_0,v_1]$ is
in the orbit of $e_1$, and $w\neq v_1$ cannot be a vertex.

This proves that $f\m(w')\subset G_v.\rond{e}_1$.
In particular, since there are no inversions on $T$, $G_{w'}\subset G_v$,
therefore $G_{w'}\subset G_{v'}$. Since $T'$ is reduced, this means that $w'$ is in the same orbit as $v'$, 
which contradicts the fact that $f\m(w')$ does not contain any vertex.

It follows that $f$ maps edges to edges, so $f$ is an isomorphism. 
In particular, $f$ induces a bijection from $V(T)$ onto $V(T')$.

The uniqueness of $f$ will then follow from the following fact:
given $v\in V(T)$, there is at most one vertex $v'\in V(T')$ such that $G_v\subset G_{v'}$.
As a matter of fact, consider $w'\in V(T')$ with $G_v\subset G_{w'}$ and let $w$ be a preimage of $w'$ in $V(T)$.
Then $\langle G_w,G_v\rangle \subset G_{w'}$ so $\langle G_w,G_v\rangle$ is elliptic in $T$.
Therefore, $v=w$ since vertex stabilizers of $T$ fix no  more than one vertex in $T$, and $v'=w'$.
\end{proof}

\def\cprime{$'$}

\end{document}